\documentclass{amsart}
\usepackage{amssymb,amsmath}

\usepackage[french]{babel}

\newtheorem{theorem}{Th\'eor\`eme}[section]
\newtheorem{lemma}[theorem]{Lemme}
\newtheorem{proposition}[theorem]{Proposition}
\newtheorem{corollary}[theorem]{Corollaire}

\theoremstyle{definition}
\newtheorem{definition}[theorem]{D\'efinition}
\newtheorem{example}[theorem]{Exemple}

\theoremstyle{remark}
\newtheorem{remark}[theorem]{Remarque}

\theoremstyle{conjecture}

\numberwithin{equation}{section}

\begin{document}

\title{Algebre Absolue} 

\author{{\bf Paul Lescot\rm}}

\address{Laboratoire de Math\'ematiques Rapha\"el Salem \\
UMR 6085 CNRS \\ 
Universit\'e de Rouen \\
Technop\^ole du Madrillet \\
Avenue de l'Universit\'e, B.P. 12 \\
76801 Saint-Etienne-du-Rouvray (FRANCE) \\ 
T\'el. 00 33 (0)2 32 95 52 24 \\
Fax 00 33 (0)2 32 95 52 86 \\
Paul.Lescot@univ-rouen.fr \\
}

\date{04 Octobre 2009}

\setcounter{section}{0}

\begin{abstract}Nous exposons la th\'eorie de Zhu concernant un analogue formel du 
corps ${\mathbf F}_{p}$, \lq\lq pour $p=1$\rq\rq, et la comparons \`a celle de Deitmar.

\bigskip

We give an exposition of Zhu's theory concerning a formal analogue of the field
${\mathbf F}_{p}$, \lq\lq for $p=1$\rq\rq, and then compare it to Deitmar's.

\end{abstract}

\maketitle

\newpage

\section{Introduction}

Il a \'et\'e, ces derni\`eres ann\'ees, propos\'e de nombreuses th\'eories du
\lq\lq corps \`a un \'el\'ement\rq\rq. Dans celle de Deitmar (\cite{4},\cite{5}), les objets de base sont les spectres de 
mono\"\i des (commutatifs, unitaires) et les sch\'emas sont obtenus par
recollement de tels objets. Le mono\"\i de trivial $F_{1}=\{1\}$ est donc l'objet final de 
la cat\'egorie de Deitmar, et son spectre $Spec(F_{1})$ l'objet initial de la cat\'egorie des $F_{1}$--sch\'emas.

Zhu (\cite{7}) a propos\'e une autre approche : $B_{1}=\{0,1\}$ d\'esigne pour lui un ensemble \`a \it deux \rm
\'el\'ements, muni de la multiplication habituelle, et d'une addition l\'eg\`erement modifi\'ee : $1+1=1$. Dans les paragraphes 2 et 3, nous d\'eveloppons, suivant Zhu mais avec des d\'emonstrations, l'alg\`ebre lin\'eaire sur $B_{1}$ en restant aussi pr\`es que possible des d\'efinitions classiques. Il s'av\`ere que la cat\'egorie des $B_{1}$--modules de type fini est beaucoup plus complexe que celle des espaces vectoriels de dimension finie sur un corps : elle est en effet \'equivalente \`a la cat\'egorie des treillis finis non vides. 

L'analogie formelle entre le groupe sym\'etrique
$\Sigma_{n}$ et le \lq\lq groupe lin\'eaire de rang $n$ sur un corps de caract\'eristique $1$\rq\rq est bien connue des sp\'ecialistes de th\'eorie des repr\'esentations ; il se trouve qu'en un sens naturel, on a bien $GL_{n}(B_{1})\simeq\Sigma_{n}$ 
(Th\-eor\`eme 3.7). Le traitement naturel (\S 4) de l'alg\`ebre commutative et de la g\'eom\'etrie alg\'ebrique sur $B_{1}$, par analogie avec l'\'etude des anneaux de polyn\^omes sur un corps, nous am\`ene \`a conclure que $B_{1}$ est \lq\lq alg\'ebriquement clos\rq\rq ,
et que, pour chaque $n\geq 1$, $MaxSpec(B_{1}[x_{1},...,x_{n}])$ est de cardinal $2^n$  (Th\'eor\`eme 4.10), ce qui constitue, dans notre cadre, un analogue du \it Nullstellensatz \rm de Hilbert.

Dans le cadre de la th\'eorie de Deitmar,  toute structure additive dispara\^\i t ;
nous faisons voir que l'on peut, au sens de la th\'eorie des cat\'egories,  plonger cette th\'eorie dans celle de Zhu, et recr\'eer ainsi une certaine
structure additive (idempotente) sur les \lq\lq anneaux de fonctions\rq\rq des objets g\'eom\'etriques consid\'er\'es.
Ce point de vue m\`ene \`a des descriptions qui nous semblent intuitivement tr\`es satisfaisantes (cf. les exemples (5.4) \`a (5.6)).

La pr\'epublication r\'ecente de Connes et Consani(\cite{2}) est post\'erieure \`a la soumission de la premi\`ere version
de cet article. Son point de vue semble, \`a certains \'egards, analogue au n\^otre. 

Les paragraphes 2 \`a 4 du texte  sont issus d'un expos\'e au Groupe de Travail Interuniversitaire en Alg\`ebre, en date du 15 Janvier 2001 ; je remercie Jacques Alev, Dominique Castella, Fran\c cois Dumas et Laurent Rigal
pour leurs commentaires \`a cette occasion. Une version pr\'eliminaire du texte complet est parue sous forme de deux pr\'epublications de l'I.H.E.S. (M/06/61 et M/06/63) \`a l'automne 2006 ; celles-ci n'auraient pas vu le jour sans l'aide et la disponibilit\'e constantes de C\'ecile Cheikhchoukh.

    J'ai pu exposer les r\'esultats de ce travail \`a l'I.H.E.S. lors de la conf\'erence \lq\lq G\'eom\'etrie Alg\'ebrique sur le Corps \`a un El\'ement\rq\rq, le 29 Mars 2007; ce m'est un agr\'eable devoir que de t\'emoigner ma gratitude \`a Christophe Soul\'e pour son invitation, et \`a plusieurs des auditeurs pour leurs commentaires et leurs encouragements, notamment Christophe Breuil, Xavier Caruso et Ofer Gabber (cf. la preuve du Th\'eor\`eme 4.8), Anton Deitmar, et Nikolai Dourov (sur la suggestion duquel j'ai chang\'e la notation de Zhu \lq\lq $F_{1}$\rq\rq en \lq\lq $B_{1}$\rq\rq afin de pr\'evenir toute confusion).

   Je remercie \'egalement le rapporteur pour de nombreuses remarques constructives, dont l'une est \`a
l'origine du contre--exemple d\'ecrit dans la Remarque 3.5, et pour avoir attir\'e mon attention sur
la pr\'epublication \cite{2}.

   Les notations sont enti\`erement standard ; si $A$ est un mono\"\i de (not\'e multiplicativement),
nous noterons $A^{*}$ l'ensemble de ses \'el\'ements inversibles.
Si $E$ est un ensemble, ${\mathcal P}_{f}(E)$ d\'enotera l'ensemble de ses parties finies, et 
\begin{eqnarray}
j_{E}
&:& E \rightarrow {\mathcal P}_{f}(E) \nonumber \\
&&x\mapsto \{x\} \nonumber 
\end{eqnarray}
l'injection canonique.

Par ${\mathcal D}$ nous entendrons la cat\'egorie des $F_{1}$--anneaux au sens de Deitmar,
c'est--\`a--dire (\cite{4}, p.88) la cat\'egorie des mono\"\i des commutatifs.
Si $A\in{\mathcal D}$ est un mono\"\i de commutatif et $B$ un sous--mono\"\i de de $A$, nous dirons que $A$ est \it entier \rm sur $B$ si, pour chaque $a\in A$, il existe un entier $n\geq 1$ tel que $a^{n}\in B$.
\section{D\'efinition de $B_{1}$ et premi\`eres propri\'t\'es}

\begin{definition}On notera $B_{1}$ l'ensemble $\{0,1\}$ muni des lois de composition internes
$+$ et $.$ donn\'ees par :
$$
0+0=0 \,\,\, ,
$$

$$
0+1=1+0=1+1=1 \,\,\, ,
$$

$$
0.0=0.1=1.0=0 \,\,\, ,
$$
et
$$
1.1=1 \,\,\, .
$$

\end{definition}

\begin{remark}Il est visible que $B_{1}$ satisfait \`a tous les axiomes des corps commutatifs,
except\'e \`a celui qui affirme l'existence de sym\'etriques pour l'addition.
\end{remark}

\begin{definition}On appelle $B_{1}$-module la donn\'ee d'un mono\"\i de commutatif $M$
d'\'el\'ement neutre $0$ et d'une $B_{1}$-loi externe sur $M$ (c'est-\`a-dire d'une application 
$$
(\lambda , x) \mapsto \lambda x
$$
de $B_{1} \times M$ dans $M$), ayant les propri\'et\'es usuelles, \it i.e. \rm :

$$
\forall (\lambda, \mu , x)\in B_{1}\times B_{1}\times M \,\,\,\, 
 (\lambda+\mu)x=\lambda x+\mu x \,\,\,\, , \leqno(1)
$$

$$
\forall (\lambda,x,y)\in B_{1}\times M\times M \,\,\, \lambda (x+y)=\lambda x+\lambda y 
\,\,\,\, , \leqno(2)
$$

$$
\forall x\in M \,\,\, 1x=x \,\,\,\, , \leqno(3)
$$

$$
\forall x\in M \,\,\, 0x=0 \,\,\,\, . \leqno(4)
$$

\end{definition} 

\begin{definition}Un ensemble ordonn\'e point\'e $(E,\leq , \, 0)$ est dit d\'ecent 
s'il poss\`ede un
(et n\'ecessairement un seul) plus petit \'el\'ement $0$, et si en outre  deux \'el\'ements 
quelconques de $E$
poss\`edent une borne sup\'erieure.
\end{definition}

\begin{theorem}La cat\'egorie $\mathcal Z$ des $B_{1}$--modules s'identifie canoniquement \`a 
la cat\'egorie des ensembles ordonn\'es d\'ecents.
\end{theorem}
\begin{proof}Soit $M$ un $B_{1}$--module ; pour $(a,b)\in M^{2}$, d\'efinissons :
$$
a\leq b \equiv a+b=b \,\, .
$$

Alors, pour tout $a\in M$ :
$$\leqno(5)$$
\begin{eqnarray}
a+a 
&=&1a+1a \nonumber \\
&=&(1+1) a \nonumber \\
&=&1a \nonumber \\
&=&a \,\,,  \nonumber 
\end{eqnarray}
soit $a\leq a$.
En outre, de $a\leq b$ et $b\leq a$, il suit :
$$
a+b=b \,\,\, \text{et} \,\,\, b+a=a \,\, ,
$$
d'o\`u $a=b+a=a+b=b$.

De plus, si $a\leq b$ et $b\leq c$, il vient :

\begin{eqnarray}
a+c
&=&a+(b+c) \nonumber \\
&=&(a+b)+c \nonumber \\
&=&b+c \nonumber \\
&=&c \,\, ,\nonumber 
\end{eqnarray}

soit $a\leq c$. $\leq$ est donc une relation d'ordre sur $M$ ; de plus, pour 
chaque $a\in M$:

$$
0+a=a \,\,\,  \,\, , \leqno(6)
$$
soit $0\leq a$ ; $(M,\leq)$ poss\`ede donc un plus petit \'el\'ement : $0$.

Soient $a\in M$ et $b\in M$ ; il est facile de voir que :

\begin{eqnarray}
a+(a+b)
&=&(a+a)+b \nonumber \\
&=&a+b \,\,\, (\text{d'apr\`es (5)}) \,\, , \nonumber 
\end{eqnarray}

soit $a\leq a+b$ ; de m\^eme $b\leq a+b$.

De plus, de $a\leq c$ et $b\leq c$ suivent $a+c=c$ et $b+c=c$, d'o\`u :

\begin{eqnarray}
(a+b)+c
&=&a+(b+c) \nonumber \\
&=&a+c \nonumber \\
&=&c \,\,\, , \nonumber
\end{eqnarray}
soit $a+b\leq c$ ; $a$ et $b$ poss\`edent donc une borne sup\'erieure :
$a\vee b=a+b$. On a bien montr\'e que $(M,\leq , 0)$ \'etait un ensemble ordonn\'e 
d\'ecent.

R\'eciproquement, soit $(E,\leq , 0)$ un ensemble ordonn\'e d\'ecent ; il est
facile de voir que l'addition et la multiplication d\'efinies par
$$
\forall (a,b)\in E^{2}\,\,\, a+b=a\vee b \,\,\, ,
$$

$$
\forall a\in E \,\,\,\,\, 0a=0 \,\,\, ,
$$
et
$$
\forall a\in E \,\,\,\,\,  1a=a \,\,\, 
$$
font de $E$ un $B_{1}$--module.

Il reste \`a d\'eterminer les \it morphismes \rm de $B_{1}$--modules. Soit $\varphi : M \rightarrow N$
un tel morphisme ; on a n\'ecessairement :
$$
\varphi(0_{M})=\varphi(0.0_{M})=0.\varphi(0_{M})=0_{N}\,\,\, ,
$$ 
et, pour $(m,m^{'})\in M^{2}$ :
$$
\varphi(m\vee_{M} m^{'})=\varphi(m+m^{'})=\varphi(m)+\varphi(m^{'})=
\varphi(m)\vee_{N}\varphi(m^{'}) \,\, .
$$
En tant qu'application entre ensembles ordonn\'es d\'ecents, $\varphi$ doit donc pr\'eserver
l'op\'eration de borne sup\'erieure (en particulier, \^etre croissante) et le plus petit
\'el\'ement.
R\'eciproquement, on v\'erifie ais\'ement qu'une application entre ensembles ordonn\'es d\'ecents
ayant ces deux propri\'et\'es constitue un morphisme pour les structures sous--jacentes de
$B_{1}$--modules.
\end{proof}

\begin{corollary}Modulo l'identification \'etablie par le Th\'eor\`eme 2.5, la cat\'egorie
$\mathcal Z_{f}$ des $B_{1}$--modules finis s'identifie \`a celle des treillis finis non vides.
\end{corollary}
\begin{proof}Soit $T$ un $B_{1}$--module fini ; par une r\'ecurrence imm\'ediate
sur le cardinal $\vert S \vert$ de $S$ on voit que toute partie (m\^eme vide) $S$ de $T$ poss\`ede une borne sup\'erieure ; en particulier,
pour $(a,b)\in T^{2}$ ,
$$
a\wedge b=\vee\{c\in T \,\, \vert \,\, c\leq a \,\,\, \text{et} \,\,\, c\leq b\}
$$
est bien d\'efini : $T$ est un treillis, et $T\neq \emptyset$ car $0\in T$.

R\'eciproquement, soit $T$ un treillis fini non vide ; il suffit de faire voir
que $T$ poss\`ede un plus petit \'el\' ement ; mais, en tant qu'ensemble ordonn\'e fini non vide,
$T$ poss\`ede un \'el\'ement minimal $m$, et on a, pour tout $x\in T$ :
$$
m\wedge x\leq m \,\,\, ,
$$
d'o\`u
$$
m\wedge x=m
$$
et
$$
m=m\wedge x\leq x \,\,\, ;
$$
$m$ est donc bien le plus petit \'el\'ement de $T$.
\end{proof}

\section{Alg\`ebre lin\'eaire sur $B_{1}$}

\begin{theorem}Soit $A$ un ensemble, et munissons l'ensemble ${\mathcal P}_{f}(A)$ de sa structure habituelle
de treillis ($C\leq B$ si et seulement si $C\subset B$) ; 
alors l'injection

\begin{eqnarray}
j_{A} : 
&& A \rightarrow {\mathcal P}_{f}(A) \nonumber \\
&& x \mapsto \{ x \} \nonumber 
\end{eqnarray}
fait de ${\mathcal P}_{f}(A)$ le $B_{1}$--module libre engendr\'e par $A$.
En particulier, la cat\'egorie des \linebreak $B_{1}$--modules libres de type fini
(i.e. finis) s'identifie canoniquement \`a celle des alg\`ebres de Boole finies.
\end{theorem}

\begin{proof}Il s'agit de faire voir que, pour tout $B_{1}$--module
$M$ et toute application $\varphi : A\rightarrow M$, il existe un unique morphisme
$$
\rho : {\mathcal P}_{f}(A)\rightarrow M
$$
tel que $\varphi=\rho \circ j_{A}$.
Pour tout $C\in {\mathcal P}_{f}(A)$, on doit avoir :

\begin{eqnarray}
\rho(C)
&=\rho(\bigcup_{x\in C}\{x\}) \nonumber \\
&=\rho(\bigcup_{x\in C}j_{A}(x)) \nonumber \\
&=\bigvee_{x\in C}\rho(j_{A}(x)) \nonumber 
\end{eqnarray}

soit :

$$
\rho(C)=\bigvee_{x\in C}\varphi(x) \,\,\,  \,\,\, , \leqno(7)
$$
d'o\`u l'unicit\'e de $\rho$.

R\'eciproquement, il est visible que $\rho$ d\'efini par (7) est un
morphisme de $B_{1}$--modules et r\'epond \`a la question.

Lorsque $A$ est fini, ${\mathcal P}_{f}(A)={\mathcal P}(A)$ est une alg\`ebre de Boole,
d'o\`u la derni\`ere assertion.
\end{proof}

Plus g\'en\'eraux que les modules libres sont les modules \it projectifs \rm,
au sens g\'en\'eral de la th\'eorie des cat\'egories : le $B_{1}$--module $M$ est
projectif si, quels que soient les $B_{1}$--modules $N_{1}$ et $N_{2}$ et les
morphismes $\varphi : M \rightarrow N_{2}$ et $\psi : N_{1} \twoheadrightarrow N_{2}$
avec $\psi$ surjectif, il existe un morphisme $\rho : M \rightarrow N_{1}$ tel que 
$\psi \circ \rho = \varphi$. Tout $B_{1}$--module libre est 
\'evidemment projectif.

\begin{definition}Soit $(E,\leq)$ un ensemble ordonn\'e ; posons
$$
{\mathcal O} (E)=\{A\subset E \, \vert \, \forall x\in A [y\leq x \Longrightarrow y\in A]\} \,\, ;
$$
alors $({\mathcal O}(E),\subset)$ est un treillis de plus petit \'el\'ement $\emptyset$, donc un
$B_{1}$--module (en fait, ${\mathcal O}(E)$ est un sous-treillis (distributif) de ${\mathcal P}(E)$).
\end{definition}

\begin{remark}Le Th\'eor\`eme suivant ne sera pas utilis\'e dans la suite de l'article.
\end{remark}

\begin{theorem}Les propri\'et\'es suivantes d'un treillis fini non vide $M$
sont \linebreak
\'equivalentes :
\begin{itemize}
\item[$\bullet$] $M$, consid\'er\'e comme $B_{1}$--module, est projectif.
\item[$\bullet$] $M$ est distributif.
\item[$\bullet$] Il existe un ensemble ordonn\'e fini $E$ tel que $M$ soit isomorphe \`a ${\mathcal O}(E)$.
\item[$\bullet$] $M$, consid\'er\'e comme $B_{1}$--module, est isomorphe \`a un sous--module d'un
\linebreak $B_{1}$--module libre.
\end{itemize}
\end{theorem}

\begin{remark}L'\'equivalence (2) $\Longleftrightarrow$ (3)
n'est autre que le cas particulier du Th\'eor\`eme de Repr\'esentation de Birkhoff 
relatif aux treillis finis : cf. par exemple \cite{1}, p.59, Theorem 3, ou \cite{3}, p.171,
Theorem 8.17.
\end{remark}
\begin{proof}
(1) $\implies$ (2) :

Soient $N_{1}={\mathcal P}(M)$, $N_{2}=M$ et 
\begin{eqnarray}
\psi :\nonumber
& {\mathcal P}(M)\,\, \rightarrow \,\, M \nonumber\\
& A \mapsto \bigvee_{x\in A}x \,\, .\nonumber
\end{eqnarray}

	Il est visible que $\psi$ est un morphisme surjectif de $B_{1}$--modules, donc il existe 
un morphisme $\rho : M \rightarrow {\mathcal P}(M)$ tel que $\psi \circ \rho=Id_{M}$.
Mais alors, pour tout $(a,b,c)\in M^{3}$:
\begin{eqnarray}
\rho (a\wedge (b \vee c))
&\leq& \rho(a)\cap \rho(b\vee c) \nonumber \\
&=&\rho(a)\cap (\rho(b) \cup \rho(c)) \nonumber \\
&=&(\rho(a)\cap \rho(b))\cup (\rho(a)\cap \rho(c)) \nonumber 
\end{eqnarray}

d'o\`u :
\begin{eqnarray}
a\wedge (b \vee c) \nonumber
&=&\psi(\rho(a\wedge(b\vee c))) \nonumber \\
&\leq& \psi((\rho(a)\cap \rho(b))\cup (\rho(a)\cap \rho(c))) \nonumber \\
&=&\psi(\rho(a)\cap \rho(b))\vee \psi(\rho(a)\cap \rho(c)) \nonumber \\
&\leq& (\psi(\rho(a))\wedge \psi(\rho(b)))\vee (\psi(\rho(a))
\wedge \psi(\rho(c))) \nonumber \\
&=&(a\wedge b)\vee (a\wedge c) \nonumber \\
&\leq& a\wedge (b\vee c) \,\,\, , \nonumber
\end{eqnarray}
donc
$$
a\wedge (b\vee c)=(a\wedge b)\vee (a\wedge c) \,\, , 
$$
c'est-\`a-dire que $\wedge$ est distributive par rapport \`a $\vee$.
Mais, comme il est bien connu (\cite{1}, Theorem 9, p.11 ; \cite{3}, Lemma 6.3, p.130), la distributivit\'e
de $\vee$ par rapport \`a $\wedge$ s'ensuit. En effet, l'on peut \'ecrire, pour
$(a,b,c)\in M^{3}$ : 

\begin{eqnarray}
(a\vee b)\wedge (a\vee c) \nonumber
&=&((a\vee b)\wedge a)\vee ((a\vee b)\wedge c) \,\,\, \text{(d'apr\`es le r\'esultat ci--dessus)} \nonumber \\
&=&a\vee (c\wedge a)\vee (c\wedge b) \,\,\, \text{(\it idem \rm)} \nonumber \\
&=&(a\vee (c\wedge a))\vee (c\wedge b) \nonumber \\
&=&a\vee (b\wedge c) \,\, . \nonumber
\end{eqnarray}

(2) $\Longrightarrow$ (3) :

Soit $E$ l'ensemble des \'el\'ements $m\neq 0$ de $M$
\it irr\'eductibles pour $\vee$ \rm, \it i.e. \rm tels que :
$$
\forall (x,y)\in M^{2} \,\,\,  x\vee y = m \Longrightarrow x=m \,\,\, \text{ou} \,\,\,y=m \,\, .
$$

De la finitude de $M$ r\'esulte que chaque \'el\'ement de $M$ est la borne sup\'erieure
d'une famille (\'eventuellement vide) d'\'el\'ements de $E$ ; dans le cas contraire, l'ensemble 
$M_{0}$ des \'el\'ements de $M$ n'ayant pas cette propri\'et\'e serait non vide, et aurait donc
un \'el\'ement minimal (pour $\leq\vert_{M_{0}}$) $a$. Par hypoth\`ese on aurait $a\neq 0$
et $a\notin E$, donc il existerait $x$ et $y$ tels que :
$$
a=x\vee y \,\, ,  x\neq a \,\, \text{et} \,\, y\neq a \,\, .
$$
Mais alors $x<a$ et $y<a$, donc $x\notin M_{0}$ et $y\notin M_{0}$, d'o\`u
$$
x=\bigvee_{b\in E_{x}}b
$$
et 
$$
y=\bigvee_{b\in E_{y}}b \,\, ,
$$
avec
$E_{x}\subset E$ et $E_{y}\subset E$. Il s'ensuivrait :

\begin{eqnarray}
a \nonumber
&=&x\vee y  \nonumber \\
&=&\bigvee_{b\in E_{x}}b \vee \bigvee_{b\in E_{y}}b \nonumber \\
&=&\bigvee_{b\in E_{x}\cup E_{y}}b\notin M_{0} \,\, , \nonumber
\end{eqnarray}
une contradiction.

On a donc :
$$
\forall m\in M \,\,\, m=\bigvee_{x\in G_{m}} x \,\, ,
$$
o\`u
$$
G_{m}=\{a\in E \,\, \vert \,\, a\leq m \} \,\, ;
$$
il est visible que $G_{m}\in{\mathcal O}(E)$. Soit alors
\begin{eqnarray}
\varphi : \nonumber
& M\rightarrow {\mathcal O}(E) \nonumber \\
& m \mapsto G_{m} \,\, ;\nonumber
\end{eqnarray}
j'affirme que $\varphi$ est un morphisme bijectif de $B_{1}$--modules.
L'injectivit\'e de $\varphi$ r\'esulte de ce que 
$$
\forall m\in M \,\,\, m=\bigvee_{x\in\varphi(m)}x \, ,
$$
la propri\'et\'e $\varphi(0)=\emptyset$ est \'evidente, et $m\leq m^{'}$ entra\^\i ne
$G_{m}\subset G_{m^{'}}$, soit

\noindent $\varphi(m)\subset \varphi(m^{'})$ ; il ne reste qu'\`a faire voir que :
$$
\varphi(m) \cup \varphi(m^{'})=\varphi(m\vee m^{'}) \,\, .
$$

L'inclusion $\varphi(m) \cup \varphi(m^{'})\subset \varphi(m\vee m^{'})$ \'etant \'evidente, il 
nous suffit d'\'etablir que:
$$
\forall x\in G_{m\vee m^{'}} \,\,\, x\in \varphi(m) \cup \varphi(m^{'})\,\, .
$$
Mais on a
\begin{eqnarray}
m\vee m^{'} \nonumber
&=&\bigvee_{a\in G_{m}}a \vee \bigvee_{a\in G_{m^{'}}}a \nonumber\\
&=&\bigvee_{a\in G_{m}\cup G_{m^{'}}}a \,\, \nonumber .
\end{eqnarray}
Soit alors $x\in G_{m\vee m^{'}}$ ; il vient :
\begin{eqnarray}
x\nonumber
&=&x\wedge (m\vee m^{'}) \nonumber\\
&=&x\wedge (\bigvee_{a\in G_{m}\cup G_{m^{'}}} a) \nonumber\\
&=&\bigvee_{a\in G_{m}\cup G_{m^{'}}}(x\wedge a) \nonumber \,\,\, .
\end{eqnarray}
Donc
$$
\exists a\in G_{m}\cup G_{m^{'}} \,\,\, x=x\wedge a \,\,\, .
$$
Mais alors $x\leq a$, d'o\`u $x\leq m$ si $a\in G_{m}$, et $x\leq m^{'}$
si $a\in G_{m^{'}}$ ; en conclusion, 
$x\in \varphi (m)$ ou $x\in \varphi(m^{'})$, et
en effet $x\in \varphi(m)\cup \varphi(m^{'})$.

Il reste maintenant \`a d\'emontrer que $\varphi(M)={\mathcal O}(E)$. Soit 
$T\in {\mathcal O}(E)$,
et soit $m=\bigvee_{t\in T} t\in M$ ; alors, pour chaque 
$t\in T$, $t\leq m$, donc $t\in \varphi(m)$ :

$$
T\subset \varphi(m) \,\, .
$$

R\'eciproquement, soit $v\in \varphi(m)$ ; on a $v\leq m$, d'o\`u :

\begin{eqnarray}
v \nonumber
&=&v \wedge m \nonumber\\
&=&v \wedge (\bigvee_{t\in T}t) \nonumber\\
&=&\bigvee_{t\in T}(v \wedge t) \nonumber 
\end{eqnarray}
donc

$$
(\exists t_{0}\in T) \,\,\,\, v=v\wedge t_{0} \,\,\, ,
$$

soit

$$
v\leq t_{0} \,\,\, ,
$$

d'o\`u (car $T\in {\mathcal O}(E)$) :

$$
v\in T \,\,\, .
$$

Il s'ensuit que $\varphi(m)\subset T$, d'o\`u

$$
T=\varphi (m)\,\,\, ;
$$

$\varphi$ est donc bel et bien surjectif.

(3) $\Longrightarrow$ (4) :

C'est \'evident vu l'existence de l'injection canonique 
$${\mathcal O}(E) \hookrightarrow {\mathcal P}(E)={\mathcal P}_{f}(E)\,\, .$$

(4) $\Longrightarrow$ (1):

On peut supposer que $M$ est un sous--$B_{1}$--module de ${\mathcal P}_{f}(E)$, pour un
certain ensemble $E$ ; en rempla\c cant \'eventuellement $E$ par $E_{1}=\bigcup_{m\in M} m$,
on peut \'egalement supposer que $E$ est fini, et que $E\in M$. Soit, pour $A\in {\mathcal P}(E)$ :
$$
{\mathcal S}(A)=\{B\in M \vert A\subset B \}.
$$
Il est clair que ${\mathcal S}(A)\neq \emptyset$ (car $E \in {\mathcal S}(A)$) ;
soit $\theta(A)=\bigcap_{B\in {\mathcal S}(A)}B$. $\theta(A)$ contient $A$ ; du fait que $M$ est un $B_{1}$--module
 \it fini \rm , donc un treillis d'apr\`es le Corollaire 2.6, r\'esulte que 
$\theta(A)\in M$ ; en particulier, $\theta(\theta(A))=\theta(A)$, \it i.e. \rm
$\theta^{2}=\theta$. Il est en outre clair que $\theta(\emptyset)=\emptyset$. 

Soient $A$ et $B$ deux \'el\'ements de ${\mathcal P}(E)$;
alors
$$
A \subset \theta(A)\subset \theta(A) \cup \theta(B) \,\, ,
$$

et de m\^eme

$$
B \subset \theta(B)\subset \theta(A) \cup \theta(B) \,\, ,
$$

soit :
$$
A\cup B \subset \theta(A)\cup \theta(B)\,\, .
$$
Mais $\theta(A) \cup \theta(B)\in M$, d'o\`u :
$$
\theta(A \cup B) \subset \theta(A) \cup \theta(B) \,\, .
$$

R\'eciproquement, si $C\in M$ et $A\cup B \subset C$ , on a $A\subset C$ et $B\subset C$,
d'o\`u $\theta(A)\subset C$ et $\theta(B)\subset C$, soit $\theta(A)\cup \theta(B) \subset C$,
d'o\`u
$$
\theta(A)\cup \theta(B)\subset \theta (A\cup B) \,\, ,
$$
et
$$
\theta(A\cup B)=\theta(A)\cup \theta(B)\,\,\,\, \,\, .
$$

Nous avons donc construit un morphisme $\theta : {\mathcal P}(E)\rightarrow M$
tel que \linebreak $\theta\vert_{M}=Id_{M}$ , c'est-\`a-dire une \it r\'etraction \rm de 
${\mathcal P}(E)$ sur $M$.
La projectivit\'e de $M$ s'ensuit alors par un raisonnement classique d'alg\`ebre universelle :
soient \linebreak $\varphi : M\rightarrow N_{2}$ et $\psi : N_{1}\twoheadrightarrow N_{2}$ surjectif deux
morphismes de $B_{1}$--modules ; alors $\varphi \circ \theta : {\mathcal P}(E)\rightarrow N_{2}$ est un 
morphisme de $B_{1}$--modules. ${\mathcal P}(E)={\mathcal P}_{f}(E)$ \'etant libre (Th\'eor\`eme 3.1), donc projectif, il existe un morphisme
$\lambda : {\mathcal P}(E)\rightarrow N_{1}$ tel que $\psi \circ \lambda = \varphi \circ \theta$.
Mais alors, en posant $\rho=\lambda\vert_{M}:M\rightarrow N_{1}$, on a :
\begin{eqnarray}
\psi \circ \rho \nonumber
&=\psi \circ \lambda\vert_{M} \nonumber \\
&=(\psi \circ \lambda)\vert_{M} \nonumber \\
&=(\varphi \circ \theta)\vert_{M} \nonumber \\
&=\varphi \circ \theta\vert_{M} \nonumber \\
&=\varphi \circ Id_{M} \nonumber \\
&=\varphi \,\,\, ; \nonumber
\end{eqnarray}
on a bien \'etabli la projectivit\'e de $M$.
\end{proof}

\begin{remark}

On peut se demander s'il est possible de caract\'eriser les treillis finis \it modulaires \rm
(lesquels forment une classe plus g\'en\'erale que celle des treillis distributifs) au moyen de leur structure de
$B_{1}$--modules, obtenant ainsi un analogue du Th\'eor\`eme 3.3. Tel ne semble pas \^etre le cas au vu de
l'exemple suivant : soit $M=\{0,a,b,c,d,e\}$ un ensemble \`a $6$ \'el\'ements, muni de la loi interne $+$ commutative, associative, idempotente,
ayant $0$ pour \'el\'ement neutre, et telle que
$$
a+b=d \,\, ,
$$
$$
c+d=e \,\, ,
$$
et
$$
b+c=c \,\, .
$$
Il est facile de v\'erifier que l'on d\'efinit ainsi sur $M$ une structure de $B_{1}$--module,
et que le treillis associ\'e (\it via \rm le Th\'eor\`eme 2.6)
est modulaire. Soit alors $N=\{0,a,c,d,e\}$ ; on voit que $N$ est un sous--$B_{1}$--module
de $M$, et que, consid\'er\'e comme treillis \it via\rm, l\`a encore, le Th\'eor\`eme 2.6,
il n'est pas modulaire (il est en fait isomorphe au treillis non--modulaire minimal $N_{5}$ :
cf. \cite{3}, (6.10), p. 134). Donc aucun \'enonc\'e portant sur un $B_{1}$--module $M$ analogue \`a la derni\`ere condition du Th\'eor\`eme 3.4
ne peut \^etre \'equivalent \`a la modularit\'e du treillis sous--jacent \`a $M$.

Dans l'exemple ci--dessus, $N$ est un sous--$B_{1}$--module de $M$, mais bien s\^ur pas un sous--treillis de ce dernier :
dans $M$
$$
c\wedge_{M} d=b
$$
et, dans $N$ :
$$
c\wedge_{N} d=0 \,\, ;
$$
aucune contradiction n'appara\^\i t donc.
\end{remark}

\begin{theorem}$GL_{n}(B_{1})\simeq \Sigma_{n}$.
\end{theorem}

\begin{proof}$GL_{n}(B_{1})$ d\'esigne par d\'efinition le groupe des 
automorphismes d'un $B_{1}$--module libre ($M$) de rang $n$. D'apr\`es le Th\'eor\`eme
3.1, on peut supposer que $M={\mathcal P_{f}}(A)={\mathcal P}(A)$ avec $\vert A \vert = n$ ; un automorphisme 
$\alpha$ de $M$ doit pr\'eserver $\emptyset$ et la relation d'inclusion, donc aussi 
les \'el\'ements minimaux de $M\setminus \{\emptyset \}$ pour l'inclusion, soit les parties 
\`a un  \'el\'ement :
$$
\forall a\in A \,\,\, \exists f(a)\in A \,\,\,  \alpha(\{a\})=\{f(a)\} \,\, .
$$
$\alpha$ \'etant injectif, l'application $f$ est injective, donc bijective, et on a, pour
tout $B\in M$ :
\begin{eqnarray}
\alpha (B)
&=&\alpha(\bigcup_{x\in B}\{x\}) \nonumber\\
&=&\bigcup_{x\in B}\alpha(\{x\}) \nonumber \\
&=&\bigcup_{x\in B}\{f(x)\} \nonumber \\
&=&\{f(x)\vert x\in B \} \nonumber \\
&=&f[B] \,\, ,\nonumber 
\end{eqnarray}
soit :

$$
\alpha(B)=f[B] \,\,\,\, \,\, . \leqno(8)
$$
R\'eciproquement, toute permutation $f$ de $A$ d\'efinit par la formule \thetag{8}
un automorphisme $\alpha$ de $M$, d'o\`u :
$$
GL_{n}(B_{1})\simeq \Sigma (A)\simeq \Sigma_{n} \,\, .
$$
\end{proof}

\section{G\'eom\'etrie alg\'ebrique sur $B_{1}$}

\begin{definition}On appelle $B_{1}$--alg\`ebre (commutative, unitaire) la donn\'ee
d'un $B_{1}$--module ${\mathcal A}$, contenant $B_{1}$, et d'une multiplication sur ${\mathcal A}$, associative,
commutative, d'\'el\'ement neutre $1$, et bilin\'eaire par rapport aux op\'erations de
$B_{1}$--module. On note $\mathcal Z_{a}$ la cat\'egorie de ces alg\`ebres.
\end{definition}

\begin{definition} On appelle congruence sur la $B_{1}$--alg\`ebre ${\mathcal A}$ 
une relation d'\'equivalence $\thicksim$ sur ${\mathcal A}$ telle que

$$
0\nsim 1
$$

et
$$
a\thicksim b \,\,\, \text{et} \,\,\, a' \thicksim b' \Longrightarrow \,\,\,
 a+a' \thicksim b+b' \,\,\, \text{et} \,\,\,  aa' \thicksim bb' \,\, .
$$
\end{definition}

Les congruences jouent dans notre th\'eorie le m\^eme r\^ole que les \'equivalences modulo
un id\'eal en alg\`ebre commutative ; en particulier, pour toute congruence 
$\thicksim$ sur ${\mathcal A}$, l'ensemble quotient ${\mathcal A}/{\thicksim}$ est muni d'une structure 
canonique de $B_{1}$--alg\`ebre.

\begin{definition}On d\'efinit sur l'ensemble des congruences sur la
$B_{1}$--alg\`ebre ${\mathcal A}$ une relation d'ordre $\geq$ par :

$$
\thicksim_{1}\geq \thicksim_{2}\iff \forall (a,b)\in {\mathcal A}^{2}\,\, 
a\thicksim_{2}b \Longrightarrow a\thicksim_{1}b \,\, .
$$
\end{definition}
Il est facile de voir que, si $\thicksim_{1}\geq \thicksim_{2}$, alors 
il existe un morphisme surjectif canonique
$$
{\mathcal A}/{\thicksim_{2}}\twoheadrightarrow {\mathcal A}/{\thicksim_{1}} \,\, .
$$

En particulier,

\begin{theorem}Si l'alg\`ebre quotient ${\mathcal A}/{\thicksim}$
est isomorphe \`a $B_{1}$, la congruence $\thicksim$ est maximale.
\end{theorem}

Il est facile de voir que la $B_{1}$--alg\`ebre libre $B_{1}[x]$
s'identifie \`a l'ensemble des sommes formelles (\'eventuellement vides) de
puissances de $x$ (en posant $x^{0}=1$). Plus g\'en\'eralement :
\begin{theorem}La $B_{1}$--alg\`ebre libre sur $A=\{x_{1},...,x_{n}\}$
s'identifie \`a l'ensemble des sommes formelles de mon\^omes 
$x_{1}^{\alpha_{1}}...x_{n}^{\alpha_{n}}$($\alpha_{i}\in \mathbf N$) muni des op\'erations
\'evidentes. Plus pr\'ecis\'ement, soit
$$
B_{1}[A]={\mathcal P}_{f}(\mathbf N^{A})
$$
l'ensemble des parties finies de $\bold N^{A}$, avec la structure naturelle de treillis,
et la multiplication d\'efinie, pour $(R,S)\in {\mathcal P}_{f}(\bold N^{A})^{2}$,  par 
$$
RS=\{a+b \,\, \vert \,\, a\in R \,\, , \,\, b\in S \}
$$
(l'addition dans $\bold N^{A}$ \'etant d\'efinie composante par composante).
Pour $a\in A$, posons $\delta_{a}=\bold 1_{\{a\}}$ ;
alors l'injection canonique
\begin{eqnarray}
i : \nonumber
& \, A \rightarrow {\mathcal P}_{f}(\bold N^{A}) \nonumber \\
& a \mapsto \{\delta_{a}\} \nonumber 
\end{eqnarray}
fait de $B_{1}[A]$ la $B_{1}$--alg\`ebre libre sur $A$.
\end{theorem}

\begin{proof}L'associativit\'e et la commutativit\'e de la multiplication sont \'evidentes,
tout comme l'existence d'un \'el\'ement neutre $U=\{0\}$ ; quant \`a la distributivit\'e, elle
suit de :
\begin{eqnarray}
R(S+T)\nonumber
&=&\{a+b \,\, \vert \,\, a\in R \,\, , \,\, b\in S+T\} \nonumber\\
&=&\{a+b \,\, \vert \,\, a\in R\, \, , \,\, b\in S\cup T\} \nonumber\\
&=&\{a+b \,\, \vert \,\, a\in R \,\, , \, \, b\in S \,\, \text{ou} \,\, b\in T\} \nonumber\\
&=&RS\cup RT \nonumber\\
&=&RS + RT \,\, . \nonumber
\end{eqnarray}

Soit maintenant $\varphi : A \rightarrow E$ une application de $A$ dans la
$B_{1}$--alg\`ebre $E$.
Il nous reste \`a montrer qu'existe un unique morphisme 
$\psi : {\mathcal P}_{f}(\bold N^{A})\rightarrow E$ tel que $\psi \circ i=\varphi$. 
Si $\psi$ est tel, on doit avoir , pour tout $F\in {\mathcal P}_{f}(\bold N^{A})$:

$$\leqno(9)$$
\begin{eqnarray}
\psi(F)\nonumber
&=&\psi(\bigcup_{x\in F}\{x\}) \nonumber\\
&=&\bigvee_{x\in F}\psi(\{x\}) \nonumber \\
&=&\bigvee_{x\in F}\psi(\{\sum_{a\in A}x(a)\delta_{a}\}) \nonumber\\
&=&\bigvee_{x\in F}\psi(\prod_{a\in A}\{\delta_{a}\}^{x(a)}) \nonumber\\
&=&\bigvee_{x\in F}\prod_{a\in A}\psi(\{\delta_{a}\})^{x(a)} \nonumber\\
&=&\bigvee_{x\in F}\prod_{a\in A}\psi(i(a))^{x(a)} \nonumber\\
&=&\bigvee_{x\in F}\prod_{a\in A}\varphi(a)^{x(a)} \nonumber \\
&=&\sum_{x\in F}\prod_{a\in A}\varphi(a)^{x(a)} \,\,\,\,   \,\, \nonumber
\end{eqnarray}
R\'eciproquement, il est tr\`es facile de voir que l'application $\psi$ d\'efinie par (9)
convient.
\end{proof}
\begin{definition}Pour $I\subset A$ et $R\in B_{1}[A]$, soit
$$
F_{I}(R)= \{r\in R \,\, \vert \,\, r(I)\subseteq\{0\}\}.
$$
\end{definition}
\begin{theorem}Pour chaque $I\subset A$, la relation $\sim_{I}$ sur $B_{1}[A]$
d\'efinie par :

\bigskip

$R\sim_{I} S$ si et seulement si ($F_{I}(R)=F_{I}(S)=\emptyset$ ou $F_{I}(R)\neq \emptyset \neq F_{I}(S)$)

\bigskip

est une congruence sur $B_{1}[A]$, et
$$
{B_{1}[A]}/{\sim_{I}}\simeq B_{1}\,\, .
$$
\end{theorem}
\begin{proof}De
$$
F_{I}(R+S)=F_{I}(R)\cup F_{I}(S)\,\, ,
$$

$$
F_{I}(RS)=F_{I}(R)F_{I}(S)\,\, ,
$$

$$
F_{I}(0)=\emptyset \,\, ,
$$
et
$$
F_{I}(1)=\{0\}=U
$$
suivent ais\'ement les propri\'et\'es qui d\'efinissent une congruence.
De plus, il est clair que $R\sim_{I} 0$ si $F_{I}(R)=\emptyset$, et que
$R\sim_{I} 1$ si $F_{I}(R)\neq \emptyset$ ; on a donc
$$
{B_{1}[A]}/{\sim_{I}}=\{\bar 0 , \bar 1\}\,\, ,
$$
d'o\`u :
$$
{B_{1}[A]}/{\sim_{I}}\simeq B_{1}.
$$
\end{proof}

En particulier, pour chaque $I\subset A$, la congruence $\sim_{I}$ sur
$B_{1}[A]$ est maximale (Th\'eor\`eme 4.4), et ${B_{1}[A]}/{\sim_{I}}\simeq B_{1}$.
R\'eciproquement, toute congruence (maximale) $\sim$ sur $B_{1}[A]$ telle que 
${B_{1}[A]}/{\sim}\simeq B_{1}$
est de la forme $\sim_{I}$ pour un $I\subset A$ (il suffit de prendre 
$$
I=\{x\in A \vert x\sim 0\}=\{x\in A \vert x\nsim 1\}).
$$

On a d'ailleurs le
\begin{theorem}Tout quotient de $B_{1}[A]$ par une congruence maximale est isomorphe \`a $B_{1}$.
\end{theorem}
\begin{proof}
Soit $\sim$ une congruence maximale sur $B_{1}[A]$. J'affirme que :

\bigskip
si $u\in B_{1}[A]$ et $v\in B_{1}[A]$ sont tels que $uv\sim 0$, alors $u\sim 0$ ou $v\sim 0$\,\,\, {(*)} ,
\bigskip

(en d'autres termes, l'alg\`ebre quotient ${B_{1}[A]}/{\sim}$ est \lq\lq int\`egre\rq\rq) .

Raisonnons par l'absurde, et soit $\mathcal R_{u}$ la relation sur $B_{1}[A]$ d\'efinie par :

$$
x \, \mathcal R_{u} \, y \equiv \exists (a,b) \in B_{1}[A]^{2} \,\,\,  x+ua \sim y+ub \,\, .
$$

Il est tr\`es facile de voir que $\mathcal R_{u}$ est compatible avec l'addition et la multiplication,
et que $x\, \sim \, y$ entra\^\i ne $x \, \mathcal R_{u} \, y$. En outre $0 \, \mathcal R_{u} \, u$, et $0\nsim u$,
donc $\sim \neq \mathcal R_{u}$. Il en r\'esulte que $\mathcal R_{u}$ n'est pas une congruence, donc que $0\, \mathcal R_{u} \, 1$, \it i.e.\rm     
il existe $(a,b)\in B_{1}[A]^ {2}$ tels que :
$$
0+ua=1+ub \,\, ,
$$
soit
$$
ua=1+ub \,\, .
$$
Mais alors

$$
(uv)a=v(ua)=v(1+ub)=v+uvb
$$

et de $uv\sim 0$ suit :
$$
0=0a\sim (uv)a=v+uvb \sim v+0b=v \,\, ,
$$

soit $v\sim 0$, une contradiction ; $(*)$ est donc bien \'etabli.

Posons $I=\{x\in A \vert x\sim 0\}$.
Soit alors $z\in B_{1}[A]$ ; d\'ecomposons $z$ en somme de mon\^omes (en les \'el\'ements de $A$) distincts :
$z = M_{1}+...+M_{k}$, et soit (pour $j\in\{1,...,k\}$) $N_{j}=_{def}\sum_{l\neq j}M_{l}$. Si $z\sim 0$, alors, pour chaque $j\in\{1,...,k\}$,

$$
0\sim z=M_{j}+N_{j}=(M_{j}+M_{j})+N_{j}=M_{j}+(M_{j}+N_{j})=M_{j}+z\sim M_{j}+0=M_{j} \,\, ,
$$
donc $M_{j}\sim 0$. Mais alors, d'apr\`es $(*)$, l'un des \'el\'ements de $A$ facteurs de $M_{j}$ est $\,\sim \,$--\'equivalent \`a 0, donc appartient \`a $I$ ;
cela valant pour chaque $j\in\{1,...,k\}$, on a $z\sim_{I} 0$.

Par ailleurs, si $z\nsim 0$ , alors $M_{j}\nsim 0$ pour au moins un $j\in\{1,...,k\}$ ;
aucun des \'el\'ements de $A$ facteurs de  $M_{j}$ n'appartient donc \`a $I$ ;
en particulier $M_{j}\sim_{I} 1$, donc $z\sim_{I} 1$.

Il en r\'esulte que $z\sim z^{'}$ entra\^\i ne $z\sim_{I} z^ {'}$, donc
que $\sim \, \leq \, \sim_{I}$, d'o\`u, au vu de la maximalit\'e de $\sim$, $\sim\, =\, \sim_{I}$.

Cette d\'emonstration est essentiellement due \`a Xavier Caruso ; j'y ai incorpor\'e quelques suggestions d'Ofer Gabber.
\end{proof}
\begin{remark}
D'apr\`es la discussion pr\'ec\'edant le Th\'eor\`eme 4.8, 
toute congruence maximale sur $B_{1}[A]$ est de la forme $\sim_{I}$ pour un $I\subset A$. Afin d'appr\'ehender
la signification de cet \'enonc\'e, consid\'erons-en l'analogue (${\mathcal E}_{K}$) sur un corps commutatif $K$ :

(${\mathcal E}_{K}$)\,\, Chaque quotient maximal de $K[x_{1},...,x_{n}]$ est isomorphe
\`a $K$, et ces quotients sont en bijection canonique avec les points de $K^{n}$.

Cet \'enonc\'e contient \`a la fois l'assertion que $K$ est alg\'ebriquement clos, et
le \it Nullstellensatz \rm . Il semble donc naturel de reformuler le Th\'eor\`eme 4.8 en le
\end{remark}

\begin{theorem}$B_{1}$ est alg\'ebriquement clos et, pour chaque $n\geq 1$, 

$Max Spec(B_{1}[x_{1},...,x_{n}])$ est de cardinal
$2^{n}$.
\end{theorem}

\begin{remark}
Les $B_{1}$--alg\`ebres \it monog\`enes \rm forment d\'ej\`a une famille
tr\`es riche. Nous nous proposons de d\'eterminer les types d'isomorphisme de $B_{1}$--alg\`ebres de cardinal $n$,
pour $n\leq 5$.
Soit donc ${\mathcal A}={B_{1}[x]}/{\sim}$ de cardinal $n$, et soit $a$ l'image de
$x\in B_{1}[x]$ dans ${\mathcal A}$ par la projection canonique.

Pour $n=2$ on a ${\mathcal A}=B_{1}$, d'o\`u
\begin{enumerate}
\item[(4.11.2.1)] 
$$
a=0
$$

ou
\item[(4.11.2.2)] 
$$
a=1 \,\, ;
$$ 
\end{enumerate}
r\'eciproquement, chacune de ces possibilit\'es d\'efinit une congruence
convenable, d'o\`u 
$$\bold{c_{2}=2\,\, . }$$

Pour $n=3$, on a n\'ecessairement $a\notin \{0,1\}$, d'o\`u
${\mathcal A}=\{0,a,1\}$.Deux cas apparaissent alors  :

\begin{enumerate} 

\item[$1^{0})$] 
$a+1=a$, soit $0<1<a$. Il suit alors $a^{2}+a=a^{2}$, d'o\`u $a^{2}\neq 1 , 0$, 
soit $a^{2}=a$, et :

\item[(4.11.3.1)] 
$$
\left\{
\aligned a+1&=a \\ 
a^{2}&=a \endaligned
\right \}
$$

\item[$2^{0})$]
$a+1=1$, soit $0<a<1$. Alors $a^{2}+a=a$, d'o\`u 
$a^{2}=0$ ou $a^{2}=a$, soit
\item[(4.11.3.2)] 
$$\left\{
\aligned a+1&=1 \\
a^{2}&=0 \endaligned \right\}
$$

ou :

\item[(4.11.3.3)]
$$\left\{
\aligned a+1&=1 \\
a^{2}&=a \endaligned \right\} .
$$

\end{enumerate}

On v\'erifie facilement que
les alg\`ebres  respectivement d\'efinies par
(4.11.3.1), (4.11.3.2) et (4.11.3.3) sont bien de cardinal 3. Il existe donc exactement trois congruences $\thicksim$
sur $B_{1}[x]$ telles que
$
{B_{1}[x]}/{\thicksim}
$
soit de cardinal 3 :
$$
\bold{c_{3}=3 \,\, .}
$$

Pour $n=4$, distinguons deux cas :
\begin{enumerate}
\item[$1^{0})$] $a^{2}\in \{0,1,a \}$. Alors $a+1\notin\{0,1,a\}$, sans quoi
$\{0,1,a\}$ serait une sous--$B_{1}$--alg\`ebre de ${\mathcal A}$ contenant $a$, et
on aurait ${\mathcal A}=\{0,1,a\}$, une contradiction. On a donc ${\mathcal A}=\{0,1,a,1+a\}$,
et $0<1<1+a$, $0<a<1+a$, et trois cas peuvent appara\^\i tre :
\item[(4.11.4.1)]
$$
a^{2}=0 \,\, ,
$$

\item[(4.11.4.2)]
$$
a^{2}=1  \,\, ,
$$

\item[(4.11.4.3)]
$$
a^{2}=a.
$$

\item[$2^{0}$)]$a^{2}\notin \{0,1,a\}$, d'o\`u ${\mathcal A}=\{0,1,a,a^{2}\}$.
Trois possibilit\'es sont alors \`a distinguer:

\item[$2^{0})\alpha)$] $a+1=1$; alors $a^{2}+a=a(a+1)=a$, d'o\`u $a^{3}+a^{2}=a(a^{2}+a)=a^{2}$, et
$0\leq a^{3} \leq a^{2}<a<1$, et encore deux \'eventualit\'es :

\item[(4.11.4.4)]
$$
\left\{ \aligned a+1&=1 \\
a^{3}&=a^{2} \endaligned \right\}\,\, ,
$$

et

\item[(4.11.4.5)]
$$
\left\{
\aligned a+1&=1 \\ 
a^{3}&=0 \endaligned \right\} \,\, .
$$

\item[$2^{0}$)$\beta$)] $a+1=a$. Alors $a^{2}+a=a(a+1)=aa=a^{2}$,
d'o\`u $a^{3}+a^{2}=a^{3}$ et $0<1<a<a^{2}\leq a^{3}$, donc $a^{2}=a^{3}$ :

\item [(4.11.4.6)]
$$
\left\{\aligned a+1&=a \\
a^{2}&=a^{3} \endaligned \right\}\,\, .
$$

\item[$2^{0}$)$\gamma$)]

\item[(4.11.4.7)] $$a^{2}=a+1\,\, .$$
\end{enumerate}
R\'eciproquement (4.11.4.1),..., (4.11.4.7) d\'efinissent chacun une alg\`ebre de cardinal 4,
d'o\`u bien :
$$
\bold{c_{4}=7 \,\, .}
$$

Pour $n=5$, distinguons \`a nouveau deux cas :
\begin{enumerate}
\item[$1^{0})$] $a+1\in \{0,1,a \}$. Alors $a^{2}\notin\{0,1,a\}$, sans quoi
$\{0,1,a\}$ serait une sous--$B_{1}$--alg\`ebre de ${\mathcal A}$ contenant $a$, et
on aurait ${\mathcal A}=\{0,1,a\}$, une contradiction
Deux cas peuvent alors se pr\'esenter :
\item[$1^{o})\alpha)$]
$$
a+1=1 \,\, .
$$
Mais alors $a^{2}+a=a(a+1)=a.1=a$, et $a^{3}+a^{2}=a^{2}(a+1)=a^{2}.1=a^{2}$,
d'o\`u $0<a^{3}<a^{2}<a<1$ ; en effet, on a n\'ecessairement $a^{3}\neq a^{2}$
et $a^{3}\neq 0$, sans quoi $\{0,1,a,a^{2}\}$ serait une sous--alg\`ebre stricte
de ${\mathcal A}$ contenant $a$. Il en r\'esulte que ${\mathcal A}=\{0, a^{3} , a^{2} , a , 1\}$ ;
du fait que $a^{4}+a^{3}=a^{3}(a+1)=a^{3}.1=a^{3}$ suit $a^{4}\leq a^{3}$
d'o\`u deux \'eventualit\'es :
\item[(4.11.5.1)]
$$
\left\{ \aligned a+1&=1 \\
a^{4}&=0 \endaligned \right\}\,\, ,
$$

et :

\item[(4.11.5.2)]
$$
\left\{ \aligned a+1&=1 \\
a^{4}&=a^{3} \endaligned \right\}\,\, .
$$

\item[$1^{o})\beta$)] $a+1=a$.

Alors $a^{2}+a=a^{2}$, $a^{3}+a^{2}=a^{3}$ , et il suit d'arguments similaires \`a ceux utilis\'es en $1^{0})\alpha)$
que $0<1<a<a^{2}<a^{3}$. Mais alors $a^{4}=a^{3}$ et :

\item[(4.11.5.3)]
$$
\left\{ \aligned a+1&=a \\
a^{4}&=a^{3} \endaligned \right\}\,\, .
$$

\item[$2^{0}$)]$a+1\notin \{0,1,a\}$.

Il s'ensuit que $a^{2}\notin \{0,1,a,a+1\}$, sans quoi
$\{0,1,a,a+1\}$ serait une sous--alg\`ebre stricte de ${\mathcal A}$ contenant $a$.
On a donc ${\mathcal A}=\{0,1,a,a+1,a^{2}\}$, et neuf possibilit\'es sont alors \`a distinguer:

\item[$2^{0})\alpha)$] $a^{2}+1=1$ et $a^{2}+a=a$ ; alors $a^{3}+a^{2}=a(a^{2}+a)=a^{2}$, et
$a^{3}+a=a(a^{2}+1)=a$, d'o\`u
$0\leq a^{3} \leq a^{2}<1$ et $0\leq a^{3} \leq a^{2} < a$ et encore deux \'eventualit\'es :

\item[(4.11.5.4)]
$$
\left\{ \aligned a^{2}+1&=1 \\
a^{2}+a&=a \\
a^{3}&=0 \endaligned \right\}\,\, ,
$$

et :

\item[(4.11.5.5)]
$$
\left\{ \aligned a^{2}+1&=1 \\
a^{2}+a&=a \\
a^{3}&=a^{2} \endaligned \right\}\,\, .
$$

\item[$2^{0})\beta)$] $a^{2}+1=1$ et $a^{2}+a=a^{2}$.

Mais alors $a+1=a+(a^{2}+1)=(a^{2}+a)+1
=a^{2}+1=1$, une contradiction.

\item[$2^{0})\gamma)$] $a^{2}+1=1$ et $a^{2}+a=a+1$.

Mais alors $a^{3}+a=a$ et $a^{3}+a^{2}=a^{2}+a=a+1$
d'o\`u $a^{3}\notin \{0,1,a+1,a^{2}\}$, et $a^{3}=a$ :

\item[(4.11.5.6)]
$$
\left\{ \aligned a^{2}+1&=1 \,\, \\
a^{2}+a&=a+1 \,\, \\
a^{3}&=a \endaligned \right\}\,\, .
$$

\item[$2^{0}) \delta )$]$a^{2}+1=a^{2}$ et $a^{2}+a=a$.

Alors $1<a^{2}<a$, d'o\`u $a+1=a$, une contradiction.

\item[$2^{0}) \epsilon )$]$a^{2}+1=a^{2}$ et $a^{2}+a=a^{2}$.

Alors il suit : $0 < 1 < a^{2}$ et $0 < a <a^{2}$, d'o\`u $a^{2}>a+1$ ; de plus $a^{3}+a^{2}=a^{3}$
d'o\`u $a^{3}\geq a^{2}$ et $a^{3}=a^{2}$ :

\item[(4.11.5.7)]
$$
\left\{ \aligned a^{2}+1&=a^{2} \\
a^{2}+a&=a^{2} \\
a^{3}&=a^{2} \endaligned \right\}\,\, .
$$

\item[$2^{0})\zeta$)]$a^{2}+1=a^{2}$ et $a^{2}+a=a+1$.

Alors $a^{3}+a=a^{3}$ et $a^{3}+a^{2}=a^{2}+a=a+1$, d'o\`u $a^{3}\notin\{0,1,a^{2}\}$,
et $a^{3}=a$ ou $a^{3}=a+1$, soit :
\item[(4.11.5.8)]
$$
\left\{ \aligned a^{2}+1&=a^{2} \\
a^{2}+a&=a+1 \\
a^{3}&=a \endaligned \right\}\,\, , \,\, \text{ou} :
$$
\item[(4.11.5.9)]
$$
\left\{ \aligned a^{2}+1&=a^{2} \\
a^{2}+a&=a+1 \\
a^{3}&=a+1 \endaligned \right\}\,\, .
$$

\item[$2^{0})\eta$)]$a^{2}+1=a+1$ et $a^{2}+a=a$.
Alors $a^{3}+a=a^{2}+a=a$ et $a^{3}+a^{2}=a^{2}$,
d'o\`u $a^{3}\notin\{1,a,a+1\}$, et $a^{3}=0$ ou $a^{3}=a^{2}$ :

\item[(4.11.5.10)]
$$
\left\{ \aligned a^{2}+1&=a+1 \\
a^{2}+a&=a \\
a^{3}&=0 \endaligned \right\}\,\, , 
$$

ou :

\item[(4.11.5.11)]
$$
\left\{ \aligned a^{2}+1&=a+1 \\
a^{2}+a&=a \\
a^{3}&=a^{2} \endaligned \right\}\,\, .
$$

\item[$2^{0})\theta )$]$a^{2}+1=a+1$ et $a^{2}+a=a^{2}$.

Alors $a^{3}+a=a^{2}+a=a^{2}$ et $a^{3}+a^{2}=a^{3}$,
d'o\`u $a^{3}\notin \{0,1,a,a+1\}$, et $a^{3}=a^{2}$, soit :
\item[(4.11.5.12)]
$$
\left\{ \aligned a^{2}+1&=a+1 \\
a^{2}+a&=a^{2} \\
a^{3}&=a^{2} \endaligned \right\}\,\, .
$$
\item[$2^{0})\iota )$]$a^{2}+1=a+1$ et $a^{2}+a=a+1$.

Alors $a^{3}+a=a^{2}+a=a+1$ et $a^{3}+a^{2}=a^{2}+a=a+1$
d'o\`u $a^{3}\notin \{0,a,a^{2}\}$, et $a^{3}=1$ ou $a^{3}=a+1$ :

\item[(4.11.5.13)]
$$
\left\{ \aligned a^{2}+1&=a+1 \\
a^{2}+a&=a+1 \\ \,\,  
a^{3}&=1 \endaligned \right\}\,\, 
$$

ou :

\item[(4.11.5.14)]
$$
\left\{ \aligned a^{2}+1&=a+1 \\
a^{2}+a&=a+1 \\
a^{3}&=a+1 \endaligned \right\}\,\, .
$$

\end{enumerate}
R\'eciproquement (4.11.5.1),..., (4.11.5.14) d\'efinissent chacun une alg\`ebre de cardinal 5,
et on a bien :
$$
\bold{c_{5}=14 \,\, .}
$$

On observe que, pour $2\leq n\leq 5$, 

$$
c_{n}=\displaystyle\frac{3}{2}n^{2}-\displaystyle\frac{13}{2}n+9 \,\, .
$$

Zhu (\cite{7}) a conjectur\'e qu'il en \'etait de m\^eme pour tout $n\geq 2$ ; il affirme
l'avoir v\'erifi\'e jusques $n=8$ inclusivement.
\end{remark}

\section{Deux foncteurs}
Comme ci--dessus, par ${\mathcal D}$ nous entendrons la cat\'egorie des $F_{1}$--anneaux au sens de Deitmar,
c'est--\`a--dire (\cite{4}, p.88) la cat\'egorie des mono\"\i des commutatifs, et par $\mathcal Z_{a}$ la cat\'egorie des $B_{1}$--alg\`ebres
au sens de la D\'efinition 4.1.
\begin{theorem}Pour $A\in {\mathcal D}$, posons ${\mathcal F}(A)={\mathcal P}_{f}(A)$, et d\'efinissons sur ${\mathcal F}(A)$ la multiplication suivante :
$$
B.C=_{def}\{xy \vert  x\in B , y\in C \}.
$$
Alors ${\mathcal F}(A)$, muni de la structure de $B_{1}$--module associ\'ee \`a sa structure
naturelle d'ensemble ordonn\'e d\'ecent (cf. Th\'eor\`eme 2.5) et de la multiplication sus-d\'efinie, constitue une $B_{1}$--alg\`ebre.
On a ${\mathcal F}(F_{1})=B_{1}$.
De plus, si $A_{1}$ et $A_{2}$ sont deux \'el\'ements de ${\mathcal D}$, et $\varphi : A_{1}\rightarrow A_{2}$ un morphisme, alors ${\mathcal F}(\varphi):{\mathcal F}(A_{1})\rightarrow {\mathcal F}(A_{2})$ d\'efini par :
$$
\forall B\in {\mathcal F}(A_{1})\,\,\,\, 
{\mathcal F}(\varphi)(B)=\{\varphi(b)\vert b\in A_{1}\}
$$
est un morphisme de ${\mathcal Z}_{a}$, et ${\mathcal F}$ d\'efinit un foncteur covariant de 
${\mathcal D}$ dans ${\mathcal Z}_{a}$.
\end{theorem}
\begin{proof}Que ${\mathcal F}(A)={\mathcal P}_{f}(A)$, muni de sa structure naturelle de treillis,
constitue un $B_{1}$--module, r\'esulte du Th\'eor\`eme 2.5.

Il est clair que la multiplication d\'efinie ci--dessus est associative, commutative, et d'\'el\'ement neutre $1_{{\mathcal F}(A)}=\{1_{A}\}$ ; sa distributivit\'e par rapport \`a l'addition est \'egalement \'evidente, l'addition dans ${\mathcal F}(A)$ n'\'etant autre  que la \it r\'eunion \rm ensembliste. De plus, pour chaque $B\in{\mathcal F}(A)$,
on a :

\begin{eqnarray}
0_{{\mathcal F}(A)}.B \nonumber
&=& \emptyset . B \nonumber\\
&=&\{ab \vert a\in \emptyset , b\in B \}\nonumber\\
&=& \emptyset \nonumber\\
&=& 0_{{\mathcal F}(A)} \,\, ; \nonumber
\end{eqnarray}
tous les axiomes de la D\'efinition 4.1 sont bien satisfaits, et ${\mathcal F}(A)$ est une $B_{1}$--alg\`ebre.
Il est clair que ${\mathcal F}(F_{1})=B_{1}$.
La validit\'e de la d\'efinition de ${\mathcal F}(\varphi)$ et la fonctorialit\'e de ${\mathcal F}$ peuvent alors \^etre v\'erifi\'ees sans probl\`eme.
\end{proof}
\begin{remark}Si $A$ est le mono\"\i de libre sur un ensemble fini $X$
(en d'autres termes, $A\simeq ({\mathbb N}^{X},+)$), alors ${\mathcal F}(A)$
s'identifie \`a la $B_{1}$--alg\`ebre libre $B_{1}[X]$ construite ci-dessus(cf. le Th\'eor\`eme 4.5, ainsi que la remarque le pr\'ec\'edant).
\end{remark}

Soit ${\mathcal G}:{\mathcal Z}_{a}\rightarrow {\mathcal D}$
le \lq\lq foncteur d'oubli\rq\rq associant \`a une $B_{1}$--alg\`ebre le mono\"\i de multiplicatif sous--jacent.

\begin{proposition}Les foncteurs ${\mathcal F}$ et ${\mathcal G}$ sont adjoints l'un de l'autre.
\end{proposition}
\begin{proof}Soient $A\in{\mathcal Z}_{a}$ et $B\in{\mathcal D}$ ; il s'agit de faire voir l'existence d'une bijection naturelle :
$$
Hom_{{\mathcal Z}_{a}}({\mathcal F}(B),A)\simeq Hom_{\mathcal D}(B,{\mathcal G}(A)) \,\, .
$$
Soit donc $\varphi : {\mathcal F}(B)\rightarrow A$ un ${\mathcal Z}_{a}$--morphisme ; pour chaque $b\in B$,
$\varphi(\{b\})$ est un \'el\'ement $\psi(b)\in A$, et il est clair que $\psi$ pr\'eserve la multiplication et l'\'el\'ement neutre (car
\begin{eqnarray}
\forall (b,b^{'})\in B^{2}\,\,\,\,  \psi(bb^{'})\nonumber
&=&\varphi(\{bb^{'}\}) \nonumber\\
&=&\varphi(\{b\}\{b^{'}\}) \nonumber\\
&=&\varphi(\{b\})\varphi(\{b^{'}\}) \nonumber\\
&=& \psi(b)\psi(b^{'}) \nonumber
\end{eqnarray}
et

$$
\psi(1_{B})=\varphi(\{1_{B}\})=\varphi(1_{{\mathcal F}(B)})=1_{A}) \,\, ,
$$
donc $\psi : B\rightarrow A={\mathcal G}(A)$ est un morphisme de ${\mathcal D}$.
D\'efinissant $\Lambda(\varphi)=\psi$, il est clair que
$$
\Lambda : Hom_{{\mathcal Z}_{a}}({\mathcal F}(B),A)\rightarrow Hom_{\mathcal D}(B,{\mathcal G}(A))
$$
est une bijection, dont la bijection inverse est d\'efinie par :

$$
\forall C\in{\mathcal F}(B)\,\,\, \Lambda^{-1}(\psi)(C)=\bigvee_{x\in C}\psi(x) \,\, .
$$
\end{proof}

\begin{example}Si $A=<x>$ est le mono\"\i de libre engendr\'e par un \'el\'ement $x$,
alors ${\mathcal F}(A)=B_{1}[x]$ est l'anneau des fonctions sur la droite affine $Spec(C_{\infty})$ de Deitmar
(\cite{4}).
\end{example}

\begin{example}Si $A=\mu_{n}$($n\geq 1$) est le groupe cyclique d'ordre $n$ ($A=<x>$ avec
$x^{n}=1$), alors ${\mathcal F}(A)$ est le quotient de $B_{1}[x]$ par la congruence
(cf. D\'efinition 4.2) engendr\'ee par la relation $x^{n}\sim 1$ : ${\mathcal F}(A)$ est de cardinal $2^{n}$, 
et appara\^\i t comme l'anneau des fonctions sur l'espace $Spec(\mu_{n})$ au sens de Deitmar
(\cite{4}).
\end{example}  

\begin{example}Si $A=<\tau , \tau^{-1}>\simeq \bold Z$ est un groupe monog\`ene infini
(engendr\'e par $\tau$), $B_{1}[A]=B_{1}[\tau, \tau^{-1}]$ est le $B_{1}$--anneau des fonctions sur le groupe multiplicatif $GL_{1}$
de Deitmar (\cite{4},p.93). On peut aussi le voir comme le quotient du $B_{1}$--anneau 
$B_{1}[x,y]$ par la congruence engendr\'ee par la relation $xy\sim 1$.
\end{example}

\begin{lemma} Pour tout $A\in {\mathcal D}$,

$$
{\mathcal F}(A)^{*}=j_{A}(A^{*})=\{\{u\}\vert u\in A^{*}\}\,\, .
$$
\end{lemma}
\begin{proof}Soit $B\in {\mathcal F}(A)$ inversible : alors il existe
$C\in {\mathcal F}(A)$ tel que 
$BC=1_{{\mathcal F}(A)}=\{1_{A}\}$.
Alors on a n\'ecessairement $B\neq \emptyset$ et $C\neq \emptyset$ ;
soient alors fix\'es $b_{0}\in B$ et $c_{0}\in C$. On a $b_{0}c_{0}\in BC$,
donc $b_{0}c_{0}=1_{A}$ : $b_{0}$ et $c_{0}$ sont inversibles.
Pour chaque $b\in B$, on a donc $bc_{0}\in BC=\{1_{A}\}$, donc
$bc_{0}=1_{A}$ et $b=c_{0}^{-1}=b_{0}$ : $B=\{b_{0}\}=j_{A}(b_{0})\in j_{A}(A^{*})$
d'o\`u $\mathcal F (A)^{*}\subseteq j_{A}(A^{*})$. L'inclusion r\'eciproque est \'evidente.
\end{proof}

\begin{corollary}Soit $\psi$ un morphisme de ${\mathcal F}(A)$ dans ${\mathcal F}(B)$;
alors 
$$
\psi(j_{A}(A^{*}))\subseteq j_{B}(B^{*})\,\, .
$$
\end{corollary}
\begin{proof}Il suffit de remarquer que $\psi$ transforme tout \'element inversible de ${\mathcal F}(A)$ en un
\'el\'ement inversible de ${\mathcal F}(B)$, et d'appliquer le Lemme 5.7.
\end{proof}

\begin{remark} Le foncteur ${\mathcal F} : {\mathcal D} \rightarrow {\mathcal Z_{a}}$
n'est pas pleinement fid\`ele : si $A=<x>$ est un mono\"\i de libre \`a un g\'en\'erateur
$x$, le morphisme de $B_{1}[x]=B_{1}[A]$ dans lui-m\^eme d\'efini par $x\rightarrow x+1$
ne provient pas d'un morphisme de mono\"\i des de $A$ dans lui-m\^eme .
\end{remark}

N\'eanmoins la situation est meilleure si l'on se restreint \`a la cat\'egorie des groupes
(ab\'eliens) :
\begin{proposition}La restriction de ${\mathcal F}$ \`a la cat\'egorie ${\mathcal A}b$ des groupes
ab\'eliens (laquelle est une sous--cat\'egorie pleine de ${\mathcal D}$) est pleinement fid\`ele.
\end{proposition}
\begin{proof}Si $A$ et $B$ sont des groupes ab\'eliens et $\psi$ un morphisme de ${\mathcal F}(A)$ dans ${\mathcal F}(B)$,
il r\'esulte de la Proposition que $\psi(j_{A}(A))\subseteq j_{B}(B)$ ; donc, pour chaque $a\in A$,
il existe $\varphi(a)\in B$ tel que $\psi(j_{A}(a))=j_{B}(\varphi(a))$. Il est maintenant clair que $\varphi : A \rightarrow B$
est un morphisme, et que $\psi={\mathcal F}(\varphi)$. L'application naturelle
$$
Hom_{{\mathcal A}b}(A,B)\rightarrow Hom_{{\mathcal Z}_{a}}({\mathcal F}(A),{\mathcal F}(B))
$$         
est donc surjective, d'o\`u le r\'esultat.
\end{proof}

\section{Une remarque}
Le Lemme suivant est l'analogue, pour des mono\"\i des, d'un r\'esultat classique de Dedekind sur les anneaux.

\begin{lemma}Soient $A\in{\mathcal D}$ un mono\"\i de commutatif, et $B\subset A$ un sous--mono\"\i de de $A$ tel que $A$ soit entier sur $B$.

Alors $A$ est un groupe si et seulement si $B$ en est un.
\end{lemma}

\begin{proof}Supposons que $A$ soit un groupe, et soit $b\in B$ ; alors il existe 
$b^{'}\in A$ tel que $bb^{'}=1$. Mais, par hypoth\`ese,
on peut trouver un entier $n\geq 1$ tel que ${b^{'}}^{n}\in B$, d'o\`u :

\begin{eqnarray}
1
&=&1^{n} \nonumber \\
&=&(bb^{'})^{n} \nonumber \\
&=&b^{n}{b^{'}}^{n} \nonumber \\
&=&b(b^{n-1}{b^{'}}^{n}). \nonumber
\end{eqnarray}
Soit $c=b^{n-1}{b^{'}}^{n}$ ; alors $c\in B$ et $bc=1$ : $b$ est donc bien inversible \it dans $B$ \rm .
Chaque \'el\'ement de $B$ y \'etant inversible, $B$ est un groupe.

R\'eciproquement, supposons que $B$ soit un groupe, et soit $a\in A$ ; par hypoth\`ese, il existe un
entier $n\geq 1$ tel que $a^{n}\in B$. $B$ \'etant un groupe, il existe $b^{'}\in B$ tel que $a^{n}b^{'}=1$ ;
mais alors :

\begin{eqnarray}
1
&=& a^{n}b^{'} \nonumber\\
&=&a(a^{n-1}b^{'}), \nonumber
\end{eqnarray}
et $a$ est inversible : $A$ est un groupe.

\end{proof}

Il s'ensuit le :
\begin{corollary} Si $B \subset A$ est une extension alg\'ebrique au sens de Deitmar 
(\cite{5}, $\S 2$), alors $A$ est un groupe si et seulement si $B$ en est un.
\end{corollary}

\bibliographystyle{amsplain}

\end{document}